\documentclass[a4paper, 11pt]{article}
\usepackage[T2A,T1]{fontenc}
\usepackage[utf8]{inputenc}
\usepackage[german, russian, english]{babel}
\usepackage{amsmath, amsthm, amssymb,amsfonts}
\usepackage{mathtools}
\usepackage{mathrsfs}
\usepackage{enumerate}
\usepackage{graphicx}
\usepackage{color}
\usepackage{hyperref}

\numberwithin{equation}{section}

\theoremstyle{definition}
\newtheorem{theorem}{Theorem}

\newtheorem{lemma}[theorem]{Lemma}
\newtheorem{remark}[theorem]{Remark}
\newtheorem{definition}[theorem]{Definition}

\newcommand{\CC}{\mathbb{C}}
\newcommand{\RR}{\mathbb{R}}
\newcommand{\NN}{\mathbb{N}}

\begin{document}

\title{Construction of Exact Control for a One-Dimensional Heat Equation with Delay}

\author{Denys Ya. Khusainov\footnote{Department of Cybernetics, Kiev National Taras Shevchenko University, Ukraine},
Michael Pokojovy\footnote{Department of Mathematics and Statistics, University of Konstanz, Germany},
Elvin I. Azizbayov\footnote{Department of Mechanics and Mathematics, Baku State University, Azerbaijan}}

\date{July 21, 2013}

\maketitle

\begin{abstract}
	We prove an exact controllability result for a one-dimensional heat equation
    with delay in both lower and highest order terms
    and nonhomogeneous Dirichlet boundary conditions.
    Moreover, we give an explicit representation of the control function
    steering the system into a given final state.
    Under certain decay properties for corresponding Fourier coefficients
    which can be interpreted as a sufficiently high Sobolev regularity of the data,
    both control function and the solution are proved to be regular
    in the classical sense both with respect to time and space variables.
\end{abstract}

\section{Introduction}
    Studying and developing mathematical models to describe various phenomena in physics, economics, ecology and population dynamics, etc.,
    are one of central problems of the modern applied mathematics (cf. \cite{EckGaKna2009}, \cite{OkuLe2001}).
    Integral and differential equations with lumped and distributed parameters proved to be a useful and efficient tool for such studies.
    Whereas evolution equations with lumped parameters have already been rather well investigated (see, e.g., \cite{Go1992}),
    there still remain a lot of open questions for the case of dynamical systems with distributed parameters
    (cp. monographs \cite{LaTri2010}, \cite{LaTri2011} and references therein).

    The scope of the present paper is a linear one-dimensional heat equation in a bounded domain
    with discrete delay in terms of both lower and highest orders.
    Recently, an abstract semigroup treatment was proposed for distributed systems with delays (viz. \cite{BaPia2001}, \cite{BaPia2005}).
    Though this rather general framework provides good analytical and control-theoretical tools
    for various delay scenarios,
    technical difficulties may arrive when applying to problems with delay in the highest order terms
    which have nevertheless been solved in \cite{BaSchn2004} for certain parabolic-type equations.

    Another important problem consists in obtaining explicit representation formulas for the solutions
    to distributed evolution equations with delay.
    We refer to \cite{AziKhu2012}, \cite{KhuIvKo2009}, \cite{KhuKu2011}, \cite{KhuPoAz2013} for details.
    Such representation formulas can then be naturally used to carefully study the solutions,
    obtain semi-analytical approximations, address controllability and optimal control problems, etc.

\section{Representation of solutions to the heat equation with delay}
    In \cite{KhuPoAz2013}, a nonhomogeneous one-dimensional heat equation with delay
    \begin{equation}
        v_{t}(x, t) = a_{1}^{2} v_{xx}(x, t) + a_{2}^{2} v_{xx}(x, t - \tau) + b_{1} v_{x}(x, t) + b_{2} v_{x}(x, t - \tau)
        + d_{1} v(x, t) + d_{2} v(x, t - \tau) + g(x, t),
        \label{NEW_EQ1}
    \end{equation}
    defined for $0 \leq x \leq l$ and $t \geq 0$ ($l > 0$), was studied.
    The coefficients for the phase derivatives were assumed to be proportional, i.e.,
    there must exist a constant $\mu \in \mathbb{R}$ such that
    $\mu = -\frac{b_{1}}{2 a_{1}^{2}} = -\frac{b_{2}}{2 a_{2}^{2}}$ holds true.
    A Dirichlet initial boundary value problem with nonhomogeneous initial
    \begin{equation}
        v(x, t) = \psi(x, t) \text{ for } 0 \leq x \leq l, -\tau \leq t \leq 0
        \label{NEW_EQ2}
    \end{equation}
    and boundary conditions
    \begin{equation}
        v(0, t) = \theta_{1}(t), \quad v(l, t) = \theta_{2}(t) \text{ for } t \geq -\tau
        \label{NEW_EQ3}
    \end{equation}
    was considered under an additional compatibility condition on the data:
    \begin{equation}
        \psi(0, t) = \theta_{1}(t), \quad \psi(l, t) = \theta_{2}(t) \text{ for } t \geq -\tau \notag
    \end{equation}
    Performing the substitution
    \begin{equation}
        v(x, t) := e^{\mu x} u(x, t) \text{ with } \mu = -\frac{b_{1}}{2 a_{1}^{2}} = -\frac{b_{2}}{2 a_{2}^{2}}, \notag
    \end{equation}
    Equation (\ref{NEW_EQ1}) was transformed to
    \begin{equation}
        u_{t}(x, t) = a_{1}^{2} u_{xx}(x, t)  + a_{2}^{2} u_{xx}(x, t - \tau) + c_{1} u(x, t) + c_{2} u(x, t) + f(x, t)
        \label{EQ1}
    \end{equation}
    with
    \begin{equation}
        c_{1} := d_{1} - \frac{b_{1}^{2}}{4 a_{1}^{2}}, \quad
        c_{2} := d_{2} - \frac{b_{2}^{2}}{4 a_{2}^{2}}, \quad
        f(x, t) := e^{-\mu x} g(x, t) \notag
    \end{equation}
    whereby the initial and boundary conditions read as
    \begin{equation}
        u(x, t) = \varphi(x, t) \text{ for } 0 \leq x \leq l, \; -\tau \leq t \leq 0, \quad
        \varphi(x, t) := e^{-\mu x} \psi(x, t)
        \label{EQ2}
    \end{equation}
    and
    \begin{equation}
        u(x, 0) = \mu_{1}(t), \quad u(l, t) := \mu_{2}(t) \text{ for } -\tau < t < 0, \quad
        \mu_{1}(t) := \theta_{1}(t), \; \mu_{2}(t) := e^{-\mu l} \theta_{2}(t),
        \label{EQ3}
    \end{equation}
    respectively.
    
\begin{figure}[h!]
        \centering
        \includegraphics[scale = 0.5]{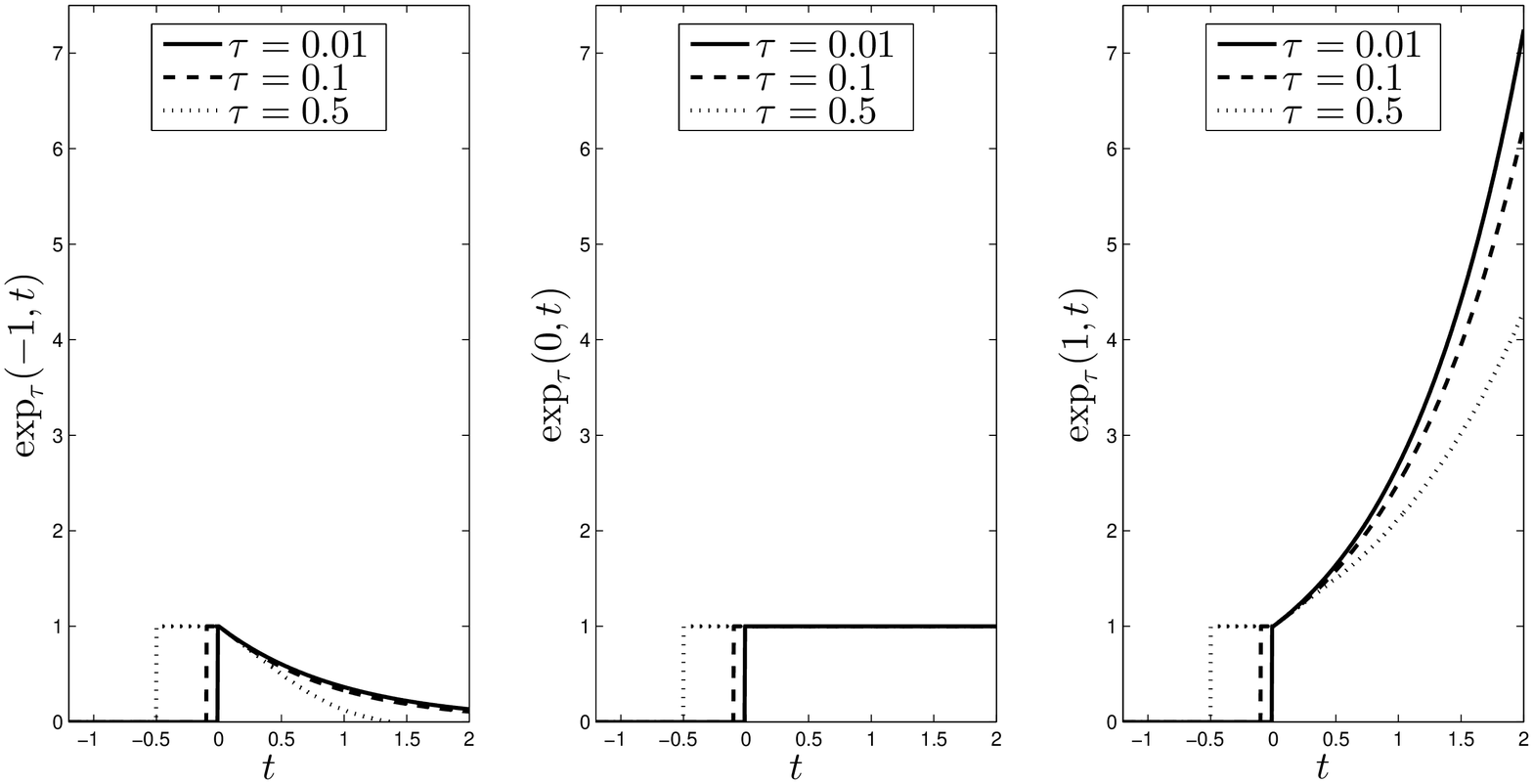}
        \caption{Delayed exponential function $\exp_{\tau}(b, \cdot)$ \label{FIGURE1}}
    \end{figure}

    Following \cite{KhuKu2011}, the delayed exponential function $\exp_{\tau}(b, \cdot)$ was introduced.
    \begin{definition}
        For $\tau > 0$, $b \in \RR$ (or $b \in \CC$), define for each $t \in \RR$:
        \begin{equation}
            \exp_{\tau}(b, t) :=
            \left\{
            \begin{array}{cc}
                0, & -\infty < t < -\tau, \\
                1, & -\tau \leq t < 0, \\
                1 + b \tfrac{t}{1!}, & 0 \leq t < \tau, \\
                1 + b \tfrac{t}{1!} + b^{2} \tfrac{(b - \tau)^{2}}{2!}, & \tau \leq t < 2\tau, \\
                \dots & \dots \\
                1 + b \tfrac{t}{1!} + \dots + b^{k} \tfrac{(t - (k - 1) \tau)^{k}}{k!}, & (k - 1) \tau \leq t < k \tau, \\
                \dots & \dots
            \end{array}\right.
            \label{EQ4}
        \end{equation}
    \end{definition}
    See Figure \ref{FIGURE1} for a plot of the delayed exponential function.

    Using the special function given in Equation (\ref{EQ4}),
    the classical solution to the initial boundary value problem (\ref{EQ1})--(\ref{EQ3})
    with delay can be represented as
    \begin{equation}
        u(x, t) = S_{1}(\varphi, \mu_{1}, \mu_{2})(x, t) + S_{2}(f, \mu_{1}, \mu_{2})(x, t) + \mu_{1}(t) + \tfrac{(\mu_{2}(t) - \mu_{1}(t)}{l} x
        \label{EQ5}
    \end{equation}
    with linear operators
    \begin{equation}
        \begin{split}
            S_{1}(\varphi, \mu_{1}, \mu_{2})(x, t) &=
            \sum_{n = 1}^{\infty} \Big(e^{L_{n}(t + \tau)} \exp_{\tau}(D_{n}, t) \Phi_{n}(-\tau) + \\
            &\phantom{=\sum_{n = 1}^{\infty} \Big(}
            \int_{-\tau}^{0} e^{L_{n}(t - s)} \exp_{\tau}(D_{n}, t - \tau - s) (\dot{\Phi}_{n}(s) - L_{n} \Phi_{n}(s)) \mathrm{d}s\Big)
            \sin(\tfrac{\pi n}{l} x), \\
            S_{2}(f, \mu_{1}, \mu_{2})(x, t) &=
            \sum_{n = 1}^{\infty} \Big(
            \int_{0}^{t} e^{L_{n}(t - s)} \exp_{\tau}(D_{n}, t - \tau - s) F_{n}(s) \mathrm{d}s\Big)
            \sin(\tfrac{\pi n}{l} x\big),
        \end{split}
        \label{EQ6}
    \end{equation}
    where
    \begin{equation}
        \begin{split}
            F_{n}(t) &= \frac{2}{l} \int_{0}^{t} f(\xi, t) \sin(\tfrac{\pi n}{l} x) \xi \mathrm{d}\xi +
            M_{n}(\mu_{1}, \mu_{2})(t), \\
            \Phi_{n}(t) &= \frac{2}{l} \int_{0}^{l} \varphi(\xi, t) \sin(\tfrac{\pi n}{l} \xi) \mathrm{d}\xi +
            m_{n}(\mu_{1}, \mu_{2})(t)
        \end{split}
        \label{EQ7}
    \end{equation}
    with
    \begin{equation}
        \begin{split}
            M_{n}(\mu_{1}, \mu_{2})(t) &= \frac{2}{l} \int_{0}^{l} \left(-\frac{\mathrm{d}}{\mathrm{d}t}
            \left(\mu_{1}(t) + \tfrac{\mu_{2}(t) - \mu_{1}(t)}{l} \xi\right) +
            \tfrac{b_{2}}{l} (\mu_{1}(t - \tau) - \mu_{2}(t - \tau))\right)
            \sin(\tfrac{\pi n}{l} \xi) \xi \mathrm{d}\xi + \\
            &\phantom{=\;\;} c_{1} m_{n}(\mu_{1}, \mu_{2})(t) + c_{2}(\mu_{1}, \mu)(t - \tau), \\
            m_{n}(\mu_{1}, \mu_{2})(t) &= \frac{2}{l} \int_{0}^{l} \left(\mu_{1}(t) + \frac{\xi}{l}(\mu_{2}(t) - \mu_{1}(t))\right)
            \sin(\tfrac{\pi n}{l} \xi) \mathrm{d}\xi
        \end{split}
        \notag
    \end{equation}
    and
    \begin{equation}
        \begin{split}
            L_{n} &= c_{1} - \left(\tfrac{\pi n}{l} a_{1}\right)^{2}, \\
            D_{n} &= \left(c_{2} - \left(\tfrac{\pi n}{l} a_{2}\right)^{2}\right)
            e^{-\left(c_{1} - \left(\tfrac{\pi n}{l} a_{1}\right)^{2}\right) \tau}.
        \end{split}
        \notag
    \end{equation}
    
    Thus, the solution to the initial boundary value problem (\ref{EQ1})--(\ref{EQ3})
    is formally given as a Fourier series in Equation (\ref{EQ6}).
    Regarding its convergence, the following result was shown in \cite{KhuPoAz2013}
    (cf. also \cite{KhuKu2011}).
    \begin{theorem}
        \label{THEOREM1}
        For $T > 0$, $m := \lceil\tfrac{T}{\tau}\rceil$ and $\alpha > 0$,
        let $F \in \mathcal{C}^{0}([0, l] \times [0, T], \RR)$, $\Phi, \partial_{t} \Phi, \partial_{tt} \Phi \in \mathcal{C}^{0}([0, l] \times [0, T], \RR)$
        be such that their Fourier coefficients $F_{n}$ and $\Phi_{n}$ satisfy
        \begin{equation}
    		\begin{split}
    			\lim\limits_{n \to \infty } n^{2m + 1 + \delta}
    			\max_{s \in [-\tau, 0]} \left[|\Phi''_{n}(s)| + n^{2} |\Phi'_{n}(s)| + n^{4} |\Phi_{n}(s)|\right] &= 0, \\
    			\lim\limits_{n \to \infty} \max_{1 \leq k \leq m} n^{2(m - k) + 1 + \delta} \max_{(k - 1) \tau \leq s \leq k\tau}
                \left[|F'_{n}(s)| + n^{2} |F_{n}(s)|\right] &= 0.
    		\end{split}
            \notag
        \end{equation}
        Under these conditions, problem (\ref{EQ1})--(\ref{EQ3}) possesses a unique classical solution
        $u \in \mathcal{C}^{0}([0, l] \times [0, T], \RR)$ with
        $\partial_{t} u, \partial_{xx} u \in \mathcal{C}^{0}([0, l] \times [0, T], \RR)$.
        Moreover, the functions $u$, $\partial_{t} u$, and $\partial_{xx} u$
        are represented by uniformly and absolutely convergent Fourier series given in (\ref{EQ5}) or obtained by a term-wise application of
        $\partial_{t}$ or $\partial_{xx}$ to (\ref{EQ5}), respectively.
    \end{theorem}

    \begin{remark}
        Using standard arguments from the elliptic theory,
        the conditions of Theorem \ref{THEOREM1} can be interpreted as a requirement for
        the data $\varphi$, $\mu_{1}$, $\mu_{2}$, $f$ to belong to
        certain Sobolev spaces (cf. \cite{Ad1975}) of functions with sufficiently many weak derivatives (s. \cite{KhuPoAz2013}).
        The larger $T$ and $\alpha$ are, the smoother the data are supposed to be.
    \end{remark}

    Representing of the solution to the initial boundary value problem (\ref{EQ1})--(\ref{EQ3}) in the form (\ref{EQ5})
    is not always convenient
    when the impact of the initial and boundary values or the inhomogeneity has to be treated separately.
    For our purposes, it is necessary to split corresponding terms into different sums.

    Expanding the first sum in (\ref{EQ5}) and performing integration by parts, we obtain
    \begin{equation}
        \begin{split}
            S_{1}&(\varphi, \mu_{1}, \mu_{2})(x, t) \\
            &= \sum_{n = 1}^{\infty}
            \left(e^{L_{n}(t + \tau)} \exp_{\tau}(D_{n}, t) \Phi_{n}(-\tau)\right)
            \sin(\tfrac{\pi n}{l} x) + \\
            &\phantom{=\,\,}
            \sum_{n = 1}^{\infty} \left(\int_{-\tau}^{0} e^{L_{n}(t - s)} \exp_{\tau}(D_{n}, t - \tau - s)
            (\dot{\Phi}_{n}(s) - L_{n} \Phi_{n}) \mathrm{d}s\right) \sin(\tfrac{\pi n}{l} x) \\
            &=
            \sum_{n = 1}^{\infty} \left(e^{L_{n} (t + \tau)} \exp_{\tau}(D_{n}, t) \Phi_{n}(-\tau)\right)
            \sin(\tfrac{\pi n}{l} x) + \\
            &\phantom{=\,\,}\sum_{n = 1}^{\infty} \left(e^{L_{n}(t - s)} \exp_{\tau}(D_{n}, t - \tau - s) \Phi_{n}(s)\Big|_{s = -\tau}^{s = 0}\right)
            \sin(\tfrac{\pi n}{l} x) - \\
            &\phantom{=\,\,}
            \sum_{n = 1}^{\infty} \left(\int_{-\tau}^{0} \left(-L_{n} e^{L_{n}(t - s)} \exp_{\tau}(D_{n}, t - \tau - s)\right) \Phi_{n}(s) \mathrm{d}s\right)
            \sin(\tfrac{\pi n}{l} x) - \\
            &\phantom{=\,\,}
            \sum_{n = 1}^{\infty} \left(\int_{-\tau}^{0} \left(-e^{L_{n}(t - s)} D_{n} \exp_{\tau}(D_{n}, t - 2\tau - s)\right) \Phi_{n}(s) \mathrm{d}s\right)
            \sin(\tfrac{\pi n}{l} x) \\
            &=
            \sum_{n = 1}^{\infty} \left(e^{L_{n} t} \exp_{\tau}(D_{n}, t - \tau) \Phi_{n}(0)
            + D_{n} \int_{-\tau}^{0} e^{L_{n}(t - s)} \exp_{\tau}(D_{n}, t - 2 \tau - s) \Phi_{n}(s) \mathrm{d}s\right)
            \sin(\tfrac{\pi n}{l} x). \notag
        \end{split}
    \end{equation}
    Plugging $\Phi_{n}$ from Equation (\ref{EQ7}), we get
    \begin{equation}
        \begin{split}
            &S_{1}(\varphi, \mu_{1}, \mu_{2})(x, t) \\
            &= \sum_{n = 1}^{\infty} \left(e^{L_{n} t} \exp_{\tau}(D_{n}, t)
            \left(\frac{2}{l} \int_{0}^{l} \varphi(\xi, 0)
            \sin(\tfrac{\pi n}{l} \xi) \mathrm{d} \xi +
            m_{n}(\mu_{1}, \mu_{2})(0)\right)\right) \sin(\tfrac{\pi n}{l} x) + \\
            &\phantom{=\,\,}
            \sum_{n = 1}^{\infty}\hspace{-0.1cm}\left(\hspace{-0.1cm}D_{n} \hspace{-0.1cm} \int_{-\tau}^{0}\hspace{-0.2cm}e^{L_{n}(t - s)} \exp_{\tau}(D_{n}, t\hspace{-0.1cm}-\hspace{-0.1cm}2 \tau\hspace{-0.1cm}-\hspace{-0.1cm}s)
            \hspace{-0.1cm}\left(\frac{2}{l}\hspace{-0.1cm}\int_{0}^{l}\hspace{-0.2cm}\varphi(\xi, s) \sin(\tfrac{\pi n}{l} x) \mathrm{d}\xi\hspace{-0.1cm}+\hspace{-0.1cm}
            m_{n}(\mu_{1}, \mu_{2})(s)\hspace{-0.1cm}\right)\hspace{-0.1cm}\mathrm{d}s\hspace{-0.1cm}\right)\hspace{-0.1cm}
            \sin(\tfrac{\pi n}{l} x).
        \end{split}
        \notag
    \end{equation}
    Expanding the sum in Equation (\ref{EQ6}) and pluggin $F_{n}$ from Equation (\ref{EQ7}) yields
    \begin{equation}
        \begin{split}
            &S_{2}(\varphi, \mu_{1}, \mu_{2})(x, t) \\
            &= \sum_{n = 1}^{\infty} \left(\int_{0}^{t} e^{L_{n}(t - s)} \exp_{\tau}(D_{n}, t - \tau - s) F_{n}(s)\right)
            \sin(\tfrac{\pi n}{l} x) \\
            &= \sum_{n = 1}^{\infty} \left(\int_{0}^{t} e^{L_{n}(t - s)} \exp_{\tau}(D_{n}, t - \tau - s)
            \left(\frac{2}{l}\hspace{-0.1cm}\int_{0}^{l}\hspace{-0.1cm}f(\xi, s) \sin(\tfrac{\pi n}{l} \xi) \mathrm{d}\xi\hspace{-0.1cm}+\hspace{-0.1cm}
            M_{n}(\mu_{1}, \mu_{2})\right)\mathrm{d}s\right)\hspace{-0.1cm} \sin(\tfrac{\pi n}{l} x).
        \end{split}
        \notag
    \end{equation}
    Thus, the solution $u$
    to the initial boundary value problem (\ref{EQ1})--(\ref{EQ3}) given in Equation (\ref{EQ5})
    can be written as follows:
    \begin{equation}
        \begin{split}
            u&(x, t) = \sum_{n = 1}^{\infty} \left(e^{L_{n} t} \exp_{\tau}(D_{n}, t - \tau) \left(\frac{2}{l} \int_{0}^{t} \varphi(\xi, 0)
            \sin(\tfrac{\pi n}{l} \xi) \right)\right) \sin(\tfrac{\pi n}{l} x) + \\
            &\phantom{=\,\,}
            \sum_{n = 1}^{\infty} \left(e^{L_{n} t} \exp_{\tau}(D_{n}, t) m_{n}(\mu_{1}, \mu_{2})(0)\right)
            \sin(\tfrac{\pi n}{l} x) + \\
            &\phantom{=\,\,}
            \sum_{n = 1}^{\infty} \left(D_{n} \int_{-\tau}^{0} e^{L_{n} (t - s)} \exp_{\tau}(D_{n}, t - 2\tau - s)
            \left(\tfrac{2}{l} \int_{0}^{l} \varphi(\xi, s) \sin(\tfrac{\pi n}{l} \xi) \mathrm{d}\xi\right)\mathrm{d}s\right)
            \sin(\tfrac{\pi n}{l} x) + \\
            &\phantom{=\,\,}
            \sum_{n = 1}^{\infty} \left(D_{n} \int_{-\tau}^{0} e^{L_{n}(t - s)} \exp_{\tau}(D_{n}, t - 2\tau - s)
            m_{n}(\mu_{1}, \mu_{2})(s) \mathrm{d}s\right) \sin(\tfrac{\pi n}{l} x) + \\
            &\phantom{=\,\,}
            \sum_{n = 1}^{\infty} \left(\int_{0}^{t} e^{L_{n}(t - s)} \exp_{\tau}(D_{n}, t - \tau - s)
            \left(\tfrac{2}{l} \int_{0}^{l} f(\xi, s) \sin(\tfrac{\pi n}{l} \xi) \mathrm{d}\xi\right) \mathrm{d}s\right)
            e^{-\tfrac{\alpha}{2} x} \sin(\tfrac{\pi n}{l} x) + \\
            &\phantom{=\,\,}
            \sum_{n = 1}^{\infty} \left(\int_{0}^{t} e^{L_{n}(t - s)} \exp_{\tau}(D_{n}, t - \tau - s)
            M_{n}(\mu_{1}, \mu_{2}) \mathrm{d}s\right) \sin(\tfrac{\pi n}{l} x) +
            \mu_{1}(t) + \tfrac{\mu_{2}(t) - \mu_{1}(t)}{l} x. \notag
        \end{split}
    \end{equation}
    Now, we collect appropriate terms in the following three operators ---
    the first one depending on the initial data:
    \begin{equation}
        \begin{split}
            \tilde{S}_{1}&(\varphi)(x, t) = \sum_{n = 1}^{\infty} \left(e^{L_{n} t} \exp_{\tau}(D_{n}, t - \tau) \left(\frac{2}{l} \int_{0}^{t} \varphi(\xi, 0)
            \sin(\tfrac{\pi n}{l} \xi) \right)\right) \sin(\tfrac{\pi n}{l} x) + \\
            &\phantom{=\,\,}
            \sum_{n = 1}^{\infty} \left(D_{n} \int_{-\tau}^{0} e^{L_{n} (t - s)} \exp_{\tau}(D_{n}, t - 2\tau - s)
            \left(\tfrac{2}{l} \int_{0}^{l} \varphi(\xi, s) \sin(\tfrac{\pi n}{l} \xi) \mathrm{d}\xi\right)\mathrm{d}s\right)
            \sin(\tfrac{\pi n}{l} x), \\
        \end{split}
        \notag
    \end{equation}
    the second one depending on the boundary data:
    \begin{equation}
        \begin{split}
            \tilde{S}_{2}&(\mu_{1}, \mu_{2})(x, t) = \sum_{n = 1}^{\infty} \left(e^{L_{n} t} \exp_{\tau}(D_{n}, t) m_{n}(\mu_{1}, \mu_{2})(0)\right)
            \sin(\tfrac{\pi n}{l} x) + \\
            &\phantom{=\,\,}
            \sum_{n = 1}^{\infty} \left(D_{n} \int_{-\tau}^{0} e^{L_{n}(t - s)} \exp_{\tau}(D_{n}, t - 2\tau - s)
            m_{n}(\mu_{1}, \mu_{2})(s) \mathrm{d}s\right) \sin(\tfrac{\pi n}{l} x) + \\
            &\phantom{=\,\,}
            \sum_{n = 1}^{\infty} \left(\int_{0}^{t} e^{L_{n}(t - s)} \exp_{\tau}(D_{n}, t - \tau - s)
            M_{n}(\mu_{1}, \mu_{2}) \mathrm{d}s\right) \sin(\tfrac{\pi n}{l} x) +
            \mu_{1}(t) + \tfrac{\mu_{2}(t) - \mu_{1}(t)}{l} x,
        \end{split}
        \notag
    \end{equation}
    and the third one depending on the inhomogeneity:
    \begin{equation}
        \begin{split}
            \tilde{S}_{2}&(f)(x, t) = \sum_{n = 1}^{\infty} \left(\int_{0}^{t} e^{L_{n}(t - s)} \exp_{\tau}(D_{n}, t - \tau - s)
            \left(\tfrac{2}{l} \int_{0}^{l} f(\xi, s) \sin(\tfrac{\pi n}{l} \xi) \mathrm{d}\xi\right) \mathrm{d}s\right)
            \sin(\tfrac{\pi n}{l} x).
        \end{split}
        \notag
    \end{equation}
    Thus, we arrive at
    \begin{equation}
        u(x, t) = \tilde{S}_{1}(\varphi)(x, t) + \tilde{S}_{2}(\mu_{1}, \mu_{2})(x, t) + \tilde{S}_{2}(f)(x, t) +
        \tfrac{\mu_{2}(t) - \mu_{1}(t)}{l} x.
        \label{EQ8}
    \end{equation}

\section{Exact controllability for the heat equation with delay}
    In this section, we consider the following exact controllability problem.
    Given an initial state $\varphi$ and boundary data $\gamma_{1}, \gamma_{2}$,
    replace $f$ with a control function $U$ such that
    the solution $u$ to (\ref{EQ1})--(\ref{EQ3}) is steered into a given final stale $\Psi$ at a prescribed time $T > 0$, i.e.,
    \begin{equation}
        u(x, T) = \Psi(x) \text{ for } 0 \leq x \leq l.
        \label{EQ9}
    \end{equation}
    Since we are interested in classical solutions,
    a compatibility condition on the boundary conditions and the end state has to be imposed:
    \begin{equation}
        \Psi(0) = \mu_{1}(T), \quad \Psi(l) = \mu_{2}(T). \notag
    \end{equation}

    As it follows from the representation formula given in Equation (\ref{EQ8}),
    Equation (\ref{EQ9}) is satisfied if and only if
    \begin{equation}
        \tilde{S}_{1}(\varphi)(x, T) + \tilde{S}_{2}(\mu_{1}, \mu_{2})(x, T) + \tilde{S}_{2}(U)(x, T) + \tfrac{\mu_{2}(T) - \mu_{1}(T)}{l} x = \Psi(x)
        \text{ for } 0 \leq x \leq l.
        \label{EQ10}
    \end{equation}

    We expand the functions $\Psi$ and $x \mapsto \tfrac{\mu_{2}(T) - \mu_{1}(T)}{l} x$
    on the interval $(0, l)$ into Fourier series with respect to the eigenfunctions of the corresponding elliptic operator.
    Equation (\ref{EQ7}) yields then
    \begin{equation}
        \begin{split}
            \Psi(x) &= \sum_{n = 1}^{\infty} \Psi_{n} \sin(\tfrac{\pi n}{l} x) \text{ with }
            \Psi_{n}(x) = \frac{2}{l} \int_{0}^{l} \Psi(\xi) \sin(\tfrac{\pi n}{l} \xi) \mathrm{d}\xi
            \text{ for } n \in \NN, \\
            \mu_{1}(t) + \tfrac{\mu_{2}(t) - \mu_{1}(t)}{l} &=
            \sum_{n = 1}^{\infty} m_{n}(\mu_{1}, \mu_{2}) \sin(\tfrac{\pi n}{l} x). \notag
        \end{split}
    \end{equation}

    We assume now $U$ to also have an expansion in Fourier series of the form:
    \begin{equation}
        U(x, t) = \sum_{n = 1}^{\infty} U_{n}(t) \sin(\tfrac{\pi n}{l} x).
        \label{EQ11}
    \end{equation}
    The operator $\tilde{S}_{3}(U)$ reads then as
    \begin{equation}
        \tilde{S}_{3}(U)(x, t) =
        \sum_{n = 1}^{\infty} \left(\int_{0}^{t} e^{L_{n} (t - s)} \exp_{\tau}(D_{n}, t - \tau - s) U_{n}(s) \mathrm{d}s\right)
        \sin(\tfrac{\pi n}{l} x). \notag
    \end{equation}
    Thus, the controllability condition rewrites as
    \begin{equation}
        \begin{split}
            \tilde{S}_{1}(\varphi)(x, T) + \tilde{S}_{2}(\mu_{1}, \mu_{2})(T, x) +
            \sum_{n = 1}^{\infty} \left(\int_{0}^{T} e^{L_{n}(t - s)} \exp_{\tau}(D_{n}, t - \tau - s) U_{n}(s) \mathrm{d}s\right)
            \sin(\tfrac{\pi n}{l} x) +& \\
            \sum_{n = 1}^{\infty} m_{n}(\mu_{1}, \mu_{2})(T) \sin(\tfrac{\pi n}{l} x) =
            \sum_{n = 1}^{\infty} \Psi_{n} \sin(\tfrac{\pi n}{l} x)&. \notag
        \end{split}
    \end{equation}
    Denote
    \begin{equation}
        \begin{split}
            s_{1n}(t) &= e^{L_{n}(t)} \exp_{\tau}(D_{n}, t) \left(\frac{2}{l} \int_{0}^{l} \varphi(\xi, 0)
            \sin(\tfrac{\pi n}{l} \xi) \mathrm{d}\xi\right) + \\
            &\phantom{=\,\,}
            D_{n} \int_{0}^{l} e^{L_{n}(t - s)} \exp_{\tau}(D_{n}, t - 2 \tau - s)
            \left(\frac{2}{l} \int_{0}^{l} \varphi(\xi, s)
            \sin(\tfrac{\pi n}{l} \xi) \mathrm{d}\xi\right) \mathrm{d} s, \\
            s_{2n}(t) &= \int_{0}^{t} e^{L_{n}(t - s)}
            \exp_{\tau}(D_{n}, t - \tau - s) M_{n}(\mu_{1}, \mu_{2})(s) \mathrm{d}s.
        \end{split}
        \notag
    \end{equation}
    There follows then from (\ref{EQ10})
    that the controllability problem for an arbitrary time $T > 0$ reduces to finding functions $u_{n}$ satisfying the following condition
    \begin{equation}
        \begin{split}
            \sum_{n = 1}^{\infty} (s_{1n}(T) + S_{2n}(T)) \sin(\tfrac{\pi n}{l} x) +
            \sum_{n = 1}^{\infty} \left(\int_{0}^{t} e^{L_{n}(t - s)} \exp_{\tau}(D_{n}, T - \tau - s) U_{n}(s) \mathrm{d}s\right)
            \sin(\tfrac{\pi n}{l} x)& + \\
            \sum_{n = 1}^{\infty} m_{n}(\mu_{1}, \mu_{2})(T)
            \sin(\tfrac{\pi n}{l} x) =
            \sum_{n = 1}^{\infty} \Psi_{n} \sin(\tfrac{\pi n}{l} x)&,
        \end{split}
        \notag
    \end{equation}
    which is in its turn equivalent to a system of countably many Fredholm integral equations of the first type:
    \begin{equation}
        s_{1n}(T) + s_{2n}(T) + \int_{0}^{T} e^{L_{n}(T - s)} \exp_{\tau}(D_{n}, T - \tau - s)
        U_{n}(s) \mathrm{d}s +
        m_{n}(\mu_{1}, \mu_{2})(T) = \Psi_{n} \text{ for } n \in \NN.
        \label{EQ12}
    \end{equation}

    \begin{lemma}
        \label{LEMMA1}
        For $\tau > 0$, $D \neq 0$, there holds for arbitrary $T > 0$
        \begin{equation}
            \int_{-\tau}^{T - \tau} \exp_{\tau}(D, s) \mathrm{d}s =
            \tfrac{1}{D} (\exp_{\tau}(D, T) - 1).
            \notag
            \label{EQ_LEMMA}
        \end{equation}
    \end{lemma}

    \begin{proof}
        There exists a unique $k \in \NN$ such that
        $(k - 2) \tau \leq T - \tau < (k - 1) \tau$. Therefore,
        \begin{equation}
            \begin{split}
                \int_{-\tau}^{T - \tau} \exp_{\tau}(D, s) \mathrm{d}s &=
                \int_{-\tau}^{0} \mathrm{d}s +
                \int_{0}^{\tau} \left(1 + D \tfrac{s}{1!}\right) \mathrm{d}s +
                \int_{\tau}^{2\tau} \left(1 + D \tfrac{s}{1!} + D^{2} \tfrac{(s - \tau)^{2}}{2!}\right) \mathrm{d}s + \\
                &\phantom{=\,\,}
                \int_{2\tau}^{3\tau} \left(1 + D \tfrac{s}{1!} + D^{2} \tfrac{(s - \tau)^{2}}{2!} + D^{3} \tfrac{(s - 2\tau)^{3}}{3!}\right) \mathrm{d}s + \dots + \\
                &\phantom{=\,\,}
                \int_{(k - 2) \tau}^{(k - 1) \tau}
                \left(1 + D \tfrac{s}{1!} + D^{2} \tfrac{(s - \tau)^{2}}{2!} + D^{3} \tfrac{(s - 2\tau)^{3}}{3!} + \dots +
                D^{k} \tfrac{(s - (k - 2)\tau)^{k-1}}{(k - 1)!}\right) \mathrm{d}s.
            \end{split}
            \notag
        \end{equation}
        Performing the integration, we obtain
        \begin{equation}
            \begin{split}
                \int_{-\tau}^{T - \tau} \exp_{\tau}(D, s) \mathrm{d}s &=
                \tfrac{s}{1!}\big|_{s = -\tau}^{s = 0} +
                \left(\tfrac{s}{1!} + D \tfrac{s^{2}}{2!}\right)\big|_{s = 0}^{s = t} +
                \left(\tfrac{s}{1!} + D \tfrac{s^{2}}{2!} + D^{2} \tfrac{(s - \tau)^{3}}{3!}\right)\big|_{s = \tau}^{s = 2\tau} + \\
                &\phantom{=\,\,}
                \left(\tfrac{s}{1!} + D \tfrac{s^{2}}{2!} + D^{2} \tfrac{(s - \tau)^{3}}{3!} + D^{3} \tfrac{(s - 2\tau)^{4}}{4!}\right)
                \big|_{s = 2\tau}^{s = 3\tau} + \dots + \\
                &\phantom{=\,\,}
                \left(\tfrac{s}{1!} + D \tfrac{s^{2}}{2!} + D^{2} \tfrac{(s - \tau)^{3}}{3!}
                + D^{3} \tfrac{(s - 2\tau)^{4}}{4!} + \dots +
                D^{k} \tfrac{(s - (k - 2) \tau)^{k}}{k!}\right)\big|_{s = (k-2)\tau}^{s = T - \tau} \\
                &= \tfrac{\tau}{1!} + \left(\tfrac{\tau}{1!} + D \tfrac{\tau^{2}}{2!}\right) +
                \left(\left(\tfrac{2\tau}{1!} + D \tfrac{(2 \tau)^{2}}{2!} + D^{2} \tfrac{\tau^{3}}{3!}\right) -
                \left(\tfrac{\tau}{1!} + D \tfrac{\tau^{2}}{2!}\right)\right) + \\
                &\phantom{=\,\,}
                \left(\left(\tfrac{3\tau}{1!} + D \tfrac{(3 \tau)^{2}}{2!} + D^{2} \tfrac{(2\tau)^{3}}{3!} + D^{3} \tfrac{\tau^{4}}{4!}\right) -
                \left(\tfrac{2\tau}{1!} + D \tfrac{(2 \tau)^{2}}{2!} + D^{2} \tfrac{\tau^{3}}{3!}\right)\right) + \dots + \\
                &\phantom{=\,\,}
                \left(\left(\tfrac{T - \tau}{1!} + D \tfrac{(T - \tau)^{2}}{2!} + D^{2} \tfrac{(T - 2\tau)^{3}}{3!}
                + D^{3} \tfrac{(T - 3\tau)^{4}}{4!} + \dots +
                D^{k} \tfrac{(T - \tau - (k - 1) \tau)^{k+1}}{(k + 1)!}\right)\right. - \\
                &\phantom{=\,\,}
                \left.\left(\tfrac{(k - 1) \tau}{1!} + D \tfrac{((k - 1) \tau)^{2}}{2!} +
                D^{2} \tfrac{((k - 2) \tau)^{3}}{3!} + D^{3} \tfrac{((k - 3) \tau)^{4}}{4!} + \dots +
                D^{k-1} \tfrac{\tau^{k}}{k!}\right)\right) \\
                &=
                \left(\tfrac{(T - \tau)}{1!} + D \tfrac{(T - \tau)^{2}}{2!} + D^{2} \tfrac{(D - 2\tau)^{3}}{3!} +
                D^{3} \tfrac{(T - 3\tau)^{4}}{4!} + \dots +
                D^{k} \tfrac{(T - \tau - (k - 1) \tau)^{k+1}}{(k + 1)!}\right) + \tfrac{\tau}{1!}.
            \end{split}
            \notag
        \end{equation}
        Thus, we can write
        \begin{equation}
            \int_{-\tau}^{T - \tau} \exp_{\tau}(D, s) \mathrm{d}s =
            \frac{1}{D} \left(1 + D \tfrac{T}{1!} + D^{2} \tfrac{(T - \tau)^{2}}{2!} +
            D^{3} \tfrac{(T - 3\tau)^{3}}{3!} + D^{4} \tfrac{(T - 3\tau)^{4}}{4!} + \dots +
            D^{k+1} \tfrac{(T - k\tau)^{k+1}}{(k + 1)!} - 1\right).
            \notag
        \end{equation}
        This completes the proof.
    \end{proof}

    Using Lemma \ref{LEMMA1},
    the integral equation (\ref{EQ12}) can be rewritten as
    \begin{equation}
        \int_{0}^{T} e^{L_{n} (T - s)} \exp_{\tau}(D_{n}, T - \tau - s) U_{n}(s) \mathrm{d}s = R_{n}(T)
        \label{EQ13}
    \end{equation}
    where
    \begin{equation}
        R_{n}(T) := \Psi_{n} - s_{1n}(T) - s_{2n}(T) - m_{n}(\mu_{1}, \mu_{2})(T). \notag
    \end{equation}
    Substituting $t := T - \tau - s$ into (\ref{EQ13}), we further obtain
    \begin{equation}
        \int_{-\tau}^{T - \tau} e^{L_{n} (\tau + t)} \exp_{\tau}(D_{n}, t)
        U_{n}(T - \tau - t) \mathrm{d}t = R_{n}(T).
        \label{EQ14}
    \end{equation}

    We look now for a solution of Equation (\ref{EQ14}) in the form
    \begin{equation}
        U_{n}(T - \tau - t) = e^{-L_{n} (\tau + t)} A_{n}(T), \notag
    \end{equation}
    where $A_{n}(T)$ are constants depending on $T$.
    Plugging this into (\ref{EQ14}) yields
    \begin{equation}
        A_{n}(T) \int_{-\tau}^{T - \tau} \exp_{\tau}(D_{n}, t) \mathrm{d}t = R_{n}(T). \notag
    \end{equation}

    Exploiting Equation (\ref{EQ_LEMMA}) from Lemma \ref{LEMMA1}, we can write
    \begin{equation}
        \tfrac{A_{n}(T)}{D_{n}} (\exp_{\tau}(D_{n}, T) - 1) = R_{n}(T). \notag
    \end{equation}
    Thus, we obtain the following Fourier coefficients for the control function
    \begin{equation}
        U_{n}(t) = e^{-L_{n}(T - t)} \tfrac{R_{n}(T) D_{n}}{\exp_{\tau}(D_{n}, T) - 1}. \notag
    \end{equation}

    Summarizing the calculations above, we have proved the following statement.
    \begin{theorem}
        Let $\varphi$, $\mu_{1}$, $\mu_{2}$ and $\Psi$ be such that the conditions of Theorem \ref{THEOREM1} are fulfilled.
        Then the control function
        \begin{equation}
            U(x, t) = \sum_{n = 1}^{\infty} U_{n}(t) \sin(\tfrac{\pi n}{l} x) \notag
        \end{equation}
        with
        \begin{equation}
            \begin{split}
                U_{n}(t) &= e^{-L_{n}(T - t)} \tfrac{R_{n}(T) D_{n}}{\exp_{\tau}(D_{n}, T) - 1}, \\
                R_{n}(T) &= \Psi_{n} - s_{1n}(T) - s_{2n}(T) - m_{n}(\mu_{1}, \mu_{2})(T), \\
                s_{1n}(t) &= e^{L_{n} t} \exp_{\tau}(D_{n}, t)
                \left(\frac{2}{l} \int_{0}^{l} \varphi(\xi, 0)
                \sin(\tfrac{\pi n}{l} \xi) \mathrm{d}\xi\right) + \\
                &\phantom{=\,\,}
                D_{n} \int_{-\tau}^{0} e^{L_{n}(t - s)} \exp_{\tau}(D_{n}, t - 2\tau - s)
                \left(\frac{2}{l} \int_{0}^{l} \varphi(\xi, s) \sin(\tfrac{\pi n}{l} \xi) \mathrm{d}\xi)\right) \mathrm{d}s, \\
                s_{2n}(t) &= \int_{0}^{t} e^{L_{n}(t - s)} \exp_{\tau}(D_{n}, t - \tau - s) M_{n}(\mu_{1}, \mu_{2})(s) \mathrm{d}s
            \end{split}
            \notag
        \end{equation}
        solves the exact controllability problem (\ref{EQ1})--(\ref{EQ3}), (\ref{EQ9}).
    \end{theorem}

    \section{Conclusions and Outlook}
        We proved an exact controllability result in the classical settings for a one-dimensional heat equation with delay.
        For practical applications, it would though be desirable to extend these results to a weak framework
        as, for example, the one described in \cite{BaSchn2004} or even beyond it.
        For $p, q \in (1, \infty)$, using the maximal $L^{p}$-regularity property (cf. \cite{Pr2002}, \cite{We2001}) for the elliptic operator in (\ref{EQ1}),
        the existence of a unique solution to (\ref{EQ1})--(\ref{EQ3})
        \begin{equation}
            u \in W^{1, p}\big((0, T), L^{q}\big((0, l)\big)\big) \cap L^{p}\big((0, T), W^{2, q}\big((0, l)\big)\big) \notag
        \end{equation}
        for the data
        \begin{equation}
            \begin{split}
                f \in L^{p}\big((0, T), L^{q}\big((0, l)\big)\big), \quad
                &\varphi \in W^{1, p}\big((-\tau, 0), L^{q}\big((0, l)\big)\big) \cap L^{p}\big((-\tau, 0), W^{2, q}\big((0, l)\big)\big), \\
                &\gamma_{1}, \gamma_{2} \in L^{p}\big((0, T), \RR\big)
            \end{split}
            \notag
        \end{equation}
        can be deduced from \cite{BaSchn2004}.
        Using this fact to verify controllability for a larger class of data and control functions
        will be a part of our further investigations.

\end{document}